\documentclass[10pt, reqno]{amsart}
\usepackage{amsmath, amsthm, amscd, amsfonts, amssymb, graphicx, color, stmaryrd}
\usepackage[all,cmtip]{xy}	
\usepackage[bookmarksnumbered, colorlinks, plainpages]{hyperref}
\usepackage[french, english]{babel}

\usepackage[onehalfspacing]{setspace}

\input{lieop.sty}

\newtheorem{Thm}{Theorem}[section]

\theoremstyle{definition}
\newtheorem{Def}[Thm]{Definition}
\newtheorem{Exa}[Thm]{Example}

\newtheorem{Rem}[Thm]{Remark}

\makeatletter
\newenvironment{abstracts}{%
  \ifx\maketitle\relax
    \ClassWarning{\@classname}{Abstract should precede
      \protect\maketitle\space in AMS document classes; reported}%
  \fi
  \global\setbox\abstractbox=\vtop \bgroup
    \normalfont\Small
    \list{}{\labelwidth\z@
      \leftmargin3pc \rightmargin\leftmargin
      \listparindent\normalparindent \itemindent\z@
      \parsep\z@ \@plus\p@
      
      \itemsep\medskipamount
    }%
}{%
  \endlist\egroup
  \ifx\@setabstract\relax \@setabstracta \fi
}

\newcommand{\abstractin}[1]{%
  \otherlanguage{#1}%
  \item[\hskip\labelsep\scshape\abstractname.]%
}
\makeatother

\begin{document}
\setcounter{page}{1}


\title[Boutet de Monvel operators on singular manifolds]{Boutet de Monvel operators on singular manifolds \\ \hfill \\ Operateurs de Boutet de Monvel pour de vari\'et\'es singuli\`eres}


\author[Karsten Bohlen]{Karsten Bohlen}

\address{$^{1}$ Leibniz University Hannover, Germany}
\email{\textcolor[rgb]{0.00,0.00,0.84}{bohlen.karsten@math.uni-hannover.de}}





\begin{abstracts}
\abstractin{english}
We construct a Boutet de Monvel calculus for general pseudodifferential boundary value problems defined on a broad class of non-compact manifolds, the class of so-called 
Lie manifolds with boundary. It is known that this class of non-compact manifolds can be used to model many classes of singular manifolds. 

\abstractin{french}
Nous construisons un calcul des type Boutet de Monvel pour des probl\`emes de valeurs au bord pseudodiff\'erentiels defin\'es sur une large classe de vari\'et\'es non-compactes, des
vari\'et\'es de Lie \`a bord. Il est bien connu que cette classe de veri\'et\'es non-compactes peut \^etre utilis\'ee pour mod\'eliser des nombreuses classes de vari\'et\'es singuli\`eres.
\end{abstracts}

\maketitle


%


\section{Introduction}

The analysis on singular manifolds has a long history, and the subject is to a large degree motivated by the study of partial differential
equations (with or without boundary conditions) and by the generalizations of index theory to the singular setting, e.g. Atiyah-Singer type index theorems.
One particular approach is based on the observation first made by A. Connes (cf. \cite{C}, section II.5) that groupoids are good models for singular spaces.
The pseudodifferential calculus on longitudinally smooth groupoids was developed by B. Monthubert, V. Nistor, A. Weinstein and P. Xu; see e.g. \cite{NWX}. 
Later a pseudodifferential calculus on a Lie manifold was constructed in \cite{ALN} via representations of pseudodifferential operators
on a Lie groupoid. This representation also yields closedness under composition.
It is important for applications in the study of partial differential equations to pose boundary conditions and to construct a parametrix for general boundary value problems.
In our case we consider the following data: a Lie manifold $(X, \V)$ with boundary $Y$ which is an embedded, transversal hypersurface $Y \subset X$ 
and which is a Lie submanifold of $X$ (cf. \cite{ALN}, \cite{AIN}). 

We will describe a general calculus with pseudodifferential boundary conditions on the Lie manifold with boundary $(X, Y, \V)$.
Special cases of our setup have been considered by Schrohe and Schulze, cf. e.g. \cite{SS}. Debord and Skandalis study Boutet de Monvel operators using deformation groupoids, \cite{DS}. 

\subsection*{Acknowledgements}

For helpful discussions and remarks I thank Magnus Goffeng, Victor Nistor, Julie Rowlett, Elmar Schrohe and Georges Skandalis. 
I thank the Deutsche Forschungsgemeinschaft (DFG) for their financial support.
%
%

\section{Boutet de Monvel's calculus}

Boutet de Monvel's calculus (e.g. \cite{BM}) was introduced in 1971. For a detailed account we refer the reader to the book \cite{G}. 
This calculus provides a convenient and general tool to study the classical boundary value problems (BVP's).
Let $X$ be a smooth compact manifold with boundary and fix smooth vector bundles $E_i \to X, \ F_i \to \partial X, \ i = 1,2$. 
Denote by $P \in \Psi_{tr}^m(M)$ a pseudodifferential operator (\cite{G}, p. 20 (1.2.4)) with transmission property (\cite{G}, p. 23, (1.2.6)) defined
on a suitable smooth neighborhood $M$ of $X$ where $P_{+}$ means $P = r^{+} P e^{+}$ the truncation. The transmission property (loc. cit.) ensures that $P_{+}$ maps functions smooth
up to the boundary to functions which are smooth up to the boundary.
Additionally, $G \colon C^{\infty}(X) \to C^{\infty}(X)$ is a singular Green operator (\cite{G}, p. 30), $K \colon C^{\infty}(\partial X) \to C^{\infty}(X)$ is a potential operator (\cite{G}, p. 29) and $T \colon C^{\infty}(X) \to C^{\infty}(\partial X)$ is a trace operator (\cite{G}, p. 27) 
We also have a pseudodifferential operator on the boundary $S \in \Psi^m(\partial X)$. 

An operator of order $m \leq 0$ and type $0$ in Boutet de Monvel's calculus is a matrix

\[
A = \begin{pmatrix} P_{+} + G & K \\
T & S \end{pmatrix} \colon \begin{matrix} C^{\infty}(X, E_1) \\ \oplus \\ C^{\infty}(\partial X, F_1) \end{matrix} \to \begin{matrix} C^{\infty}(X, E_2) \\ \oplus \\ C^{\infty}(\partial X, F_2) \end{matrix} \in \B^{m,0}(X, \partial X).
\]

%

\textbf{The calculus of Boutet de Monvel has the following \emph{features:}} 
\begin{itemize}
\item If the bundles match, i.e. if $E_1 = E_2 = E, \ F_1 = F_2 = F$ the calculus is \textbf{\emph{closed under composition}}. 

\item If $F_1 = 0, \ G = 0$ and $K, \ S$ are not present, we obtain a \textbf{\emph{classical BVP}}, e.g. the Dirichlet problem.

\item If $F_2 = 0$ and $T, \ S$ are not present, the calculus \textbf{\emph{contains inverses of classical BVP's}} whenever they exist. 
\end{itemize}

The proof that two Boutet de Monvel operators composed are again of this type is technical, see e.g. \cite{G}, chapter 2. 
Additionally, the symbolic structure of the operators involves more complicated behavior than that of merely pseudodifferential operators.

\section{Lie manifolds with boundary}

In this section we will consider the general setup for the analysis on singular and non-compact manifolds.

\begin{Exa}
On a compact manifold with boundary $M$ we introduce a Riemannian metric which models a singular structure, where the manifold with boundary is viewed as a compactification
of a non-compact manifold with cylindrical end.
Precisely, the metric is a product metric $g = g_{\partial M} + dt^2$ in a tubular neighborhood of the boundary (or the far end of
the cylinder). 
The cylindrical end is mapped to a tubular neighborhood of the boundary via the Kondratiev transform $r = e^t$ based on \cite{Kon}. 
Assume that we are given a tubular neighborhood of the form $[0, \epsilon) \times \partial M$ and let $(r, x') \in [0,\epsilon) \times \partial M$
be local coordinates.
The $b$-\emph{differential operators} take the form for $n = \dim(M)$
\begin{align}
P &= \sum_{|\alpha| \leq m} a_{\alpha}(r, x') (r \partial_r)^{\alpha_1} \partial_{x_2'}^{\alpha_2} \cdots \partial_{x_n'}^{\alpha_n} = \sum_{|\alpha| \leq m} a_{\alpha} (r \partial_r)^{\alpha_1} \partial^{\alpha'}. \tag{$*$} \label{1}
\end{align}

We observe that the vector fields, which are \emph{local generators}, in this example are $\{r \partial_r, \partial_{x_2}, \cdots, \partial_{x_n}\}$.

We consider a locally finitely generated module of vector fields $\V_b$ that has local generators as defined in our example.
These are the vector fields that are tangent to the boundary $\partial M$. 
An operator of order $m$ in the universal enveloping algebra $P \in \Diff_{\V_b}^m(M)$ is locally written as in \eqref{1}. 
Since $\V_b$ is a locally finitely generated and projective $C^{\infty}(M)$-module we obtain a vector bundle $\A_b \to M$ such that the smooth sections identify $\Gamma(\A_b) \cong \V_b$ by the Serre-Swan theorem. On $\A_b$ we have a structure of a \emph{Lie algebroid} with anchor $\varrho \colon \A_b \to TM$ (see \cite{ALN} for further details).
\label{Exa:Kontradiev}
\end{Exa}

There are sub Lie algebras of $\V_b$ constituting so-called \emph{Lie structures} which model different types of singular structures on a manifold, see also \cite{AIN}, \cite{ALN}, \cite{ALNV}.
In this general setup $M$ is a compact manifold \emph{with corners} (generalizing Example \ref{Exa:Kontradiev}) viewed as the compactification endowed with a Riemannian metric. 
The Riemannian metric is of product type in a tubular neighborhood of the singular hyperfaces. The topological structure of $M$ is such that $M$ has a finite number of embedded (intersecting) codimension one hypersurfaces. 
Open subsets of $[-1,1]^k \times \Rr^{n-k}$, where $k$ is the codimension, are needed to model manifolds with corners.
We can require the transition maps to be smooth and obtain a smooth structure on $M$. In our setup we consider such a Lie manifold $X$ with an additional hypersurface (denoted $Y$ below) which is \emph{transversal}.
Hence $Y$ is allowed to intersect the singular strata (at infinity) of $X$ as long as this intersection does not
occur in a corner (where two singular strata meet).

\begin{Def}[\cite{ALN}, Def. 1.1]
A \emph{Lie manifold} $(X, \A^{\pm})$ consists of the following data.

\emph{i)} A compact manifold with corners $X$. 

\emph{ii)} A Lie algebroid $(\A^{\pm}, \varrho_{\pm})$ with projection map $\pi_{\pm} \colon \A^{\pm} \to X$. 

\emph{iii)} The module of vector fields $\V_{\pm} = \Gamma(\A^{\pm})$ is a locally finitely generated, projective
$C^{\infty}(X)$-module. 
\end{Def}

\begin{Def}[\cite{AIN}, Def. 2.1, Def. 2.5]
A \emph{Lie manifold with boundary} $(X, Y, \A^{\pm})$ consists of the following data.

\emph{i)} A Lie manifold $(X, \A^{\pm})$. 

\emph{ii)} An embedded codimension one submanifold with corners $Y \hookrightarrow X$. 

\emph{iii)} There is a Lie algebroid $(\A_{\partial}, \varrho_{\partial})$ on $Y$ with projection map
$\pi_{\partial} \colon \A_{\partial} \to Y$ such that $\A_{\partial}$ is a Lie subalgebroid of $\A^{\pm}$, \cite{M}, Def. 4.3.14. 

\emph{iv)} The submanifold $Y$ is \emph{transversal}, i.e. $\varrho_{\pm}(\A_y) + T_y Y = T_y X, \ y \in \partial Y$. 

\emph{v)} The interior $(X_0, Y_0)$ is diffeomorphic to a smooth manifold with boundary. 
\end{Def}

\begin{Rem}
\emph{i)} In \cite{AIN} the authors define a Lie manifold with boundary $(X, Y, \V_{\pm})$ and also the \emph{double} of a given Lie manifold with boundary.
We denote this double by $M = 2X$ which is a Lie manifold $(M, \V)$.
The \emph{Lie structure} $\V$ is defined such that $\V_{\pm} = \{V_{|X_{\pm}} : V \in \V\}$.
We obtain a Lie manifold $(M, \A)$ with the Lie algebroid $(\A, \varrho)$. 

\emph{ii)} We set $\W = \Gamma(\A_{\partial})$ for the \emph{Lie structure} of $Y$. Using \emph{iii)} and \emph{iv)}
of the definition we obtain (cf. \cite{M}, pp. 164-165)
\begin{align*}
\W &= \{V \in \Gamma(Y, \A_{|Y}) : \varrho \circ V \in \Gamma(Y, TY)\} = \{V_{|Y} : V \in \V, \ V_{|Y} \ \text{tangent to} \ Y\}. 
\end{align*}

\label{Rem:double}
\end{Rem}

\section{Quantization}


In this section we describe the quantization of H\"ormander symbols defined on the conormal bundles.
We restrict ourselves to the case of trace operators. The other cases are defined analogously.

\textbf{\emph{We fix the following data:}} 
\begin{itemize}
\item A Lie manifold with boundary $(X, Y, \V)$ and the double $M = 2X$ of $X$, endowed with Lie structure $2 \V$. Fix a Lie algebroid $(\pi \colon \A \to M, \ \varrho_M)$ such that $\Gamma(\A) = 2\V$. The hypersurface $Y$ is endowed with the Lie structure $\W$ as defined in \ref{Rem:double}.
Furthermore, fix the vector bundle $(\pi_{\partial} \colon \A_{\partial} \to Y, \ \varrho_{\partial})$ with $\Gamma(\A_{\partial}) = \W$.  

\item We fix the \emph{normal bundles} $\A_{|Y} / \A_{\partial} =: \N \to Y$ as well as\footnote{We denote by $\Delta_Y$ the diagonal in $Y \times Y$ being understood as a submanifold of $Y \times M, \ M \times Y$ and $M \times M$, groupoids $\G$, $\G_{\partial}$ and spaces $\X$, $\Xop$ as defined in \cite{B}.} $\N^{\X} \Delta_Y \to Y, \ \N^{\X^t} \Delta_Y \to Y, \ \N^{\G} \Delta_Y \to Y$ which are used to quantize pseudodifferential, trace, potential and singular Green operators respectively. 
\end{itemize}

The notation is reminiscent of the underlying geometry which is described using groupoids and groupoid correspondences \cite{B}.
We will keep using this notation, though we remark that there are (non-canonical) isomorphisms
\begin{align*}
& \N^{\X} \Delta_Y \cong \A_{\partial} \times \N, \ \N^{\Xop} \Delta_{Y} \cong \N \times \A_{\partial} \ \text{and} \ \N^{\G} \Delta_Y \cong \A_{|Y} \times \N. 
\end{align*}

\begin{Rem}
\emph{i)} On the singular normal bundles we define the H\"ormander symbols spaces $S^m(\N^{\X} \Delta_Y^{\ast}) \subset C^{\infty}(\N^{\X} \Delta_Y^{\ast})$ as in \cite{HIII}, Thm 18.2.11.

\emph{ii)} Define the \emph{inverse fiberwise Fourier transform} 
\[
\Ff^{-1}(\varphi)(\zeta) = \int_{\overline{\pi}(\zeta) = \pi(\xi)} e^{i \scal{\xi}{\zeta}} \varphi(\xi) \,d\xi, \ \varphi \in S(\N^{\X} \Delta_Y^{\ast}).
\]

Here we use the notation $S(\N^{\X} \Delta_Y^{\ast})$ for the space of rapidly decreasing functions on the conormal bundle, see also \cite{S}, Chapter 1.5.

The spaces of conormal distributions are defined as $I^{m}(\N^{\X} \Delta_Y, \Delta_Y) := \Ff^{-1} S^m(\N^{\X} \Delta_Y^{\ast})$ and $I^{m}(\N^{\Xop} \Delta_Y, Y), \ I^{m}(\N^{\G} \Delta_Y, Y)$ analogously.
\label{Rem:fwise}
\end{Rem}

On a Lie manifold the injectivity radius is positive, see \cite{ALN2}, Thm. 4.14. 
Let $r$ be smaller than the injectivity radius and write $(\N^{\X} \Delta_Y)_r = \{v \in \N^{\X} \Delta_Y : \|v\| < r\}$ as well as $I_{(r)}^m(\N^{\X} \Delta_Y, \Delta_Y) = I^m((\N^{\X} \Delta_Y)_r, \Delta_Y).$

Fix the restriction $\R \colon I_{(r)}^m(\N^{\X} \Delta_Y, \Delta_Y) \to I_{(r)}^m(N^{Y_0 \times M_0} \Delta_{Y_0}, \Delta_{Y_0})$.
Additionally, denote by $\J_{tr}$ the action of a conormal distribution (its induced linear operator).
We denote by $\Psi$ the normal fibration of the inclusion $\Delta_{Y_0} \hookrightarrow Y_0 \times M_0$
such that $\Psi$ is the local diffeomorphism mapping an open neighborhood of the zero section $O_{Y_0} \subset V \subset N^{Y_0 \times M_0} \Delta_{Y_0}$
onto an open neighborhood $\Delta_{Y_0} \subset U \subset Y_0 \times M_0$ (cf. \cite{S}, Thm. 4.1.1). 
Then we have the induced map on conormal distributions $\Psi_{\ast} \colon I_{(r)}^m(N^{Y_0 \times M_0} \Delta_{Y_0}, \Delta_{Y_0}) \to I^m(Y_0 \times M_0, \Delta_{Y_0})$.
Also let $\chi \in C_c^{\infty}(\N^{\X} \Delta_Y)$ be a cutoff function which acts by multiplication $I^m(\N^{\X} \Delta_Y, \Delta_Y) \to I_{(r)}^{m}(\N^{\X} \Delta_Y, \Delta_Y)$.

\begin{Def}[Quantization]
Define $q_{T, \chi} \colon S^m(\N^{\X} \Delta_Y^{\ast}) \to \Trace^{m,0}(M, Y)$ such that for $t \in S^m(\N^{\X} \Delta_Y^{\ast})$ we have $q_{T, \chi}(t) = \J_{tr} \circ q_{\Psi, \chi}(t)$ where $q_{\Psi, \chi}(t) = \Psi_{\ast}(\R(\chi \Ff^{-1}(t)))$.
\end{Def}

From the compactness of $M$ we can associate to each vector field in $2 \V$ a \emph{global flow} $2\V \ni V \mapsto \Phi_V \colon \Rr \times M \to M$.
Then consider the diffeomorphism $\Phi(1, -) \colon M \to M$ evaluated at time $t = 1$ and fix the corresponding group actions on functions which we denote by $2\V \ni V \mapsto \varphi_V \colon C^{\infty}(M) \to C^{\infty}(M)$.

\begin{Def}
The class of $\V$-trace operators is defined as $\Trace_{2\V}^{m,0}(M, Y) := \Trace^{m,0}(M, Y) + \Trace_{2\V}^{-\infty, 0}(M, Y)$.
Here $\Trace^{m,0}(M, Y)$ consists of the extended operators from the previous definition.
The residual class is defined as follows
\begin{align*}
& \Trace_{2\V}^{-\infty, 0}(M, Y) := \mathrm{span}\{q_{\chi, T}(t) \varphi_{V_1} \cdots \varphi_{V_k} : V_j \in 2\V, \ \chi \in C_c^{\infty}(\N^{\X} \Delta_Y), \ t \in S^{-\infty}(\N^{\X} \Delta_Y^{\ast})\}.
\end{align*}
\end{Def}

We henceforth denote by $\B_{2\V}^{m,0}(M, Y)$ the class of extended Boutet de Monvel operators which consist of
matrices of operators $\begin{pmatrix} P + G & K \\ T & S \end{pmatrix}$. The components are given via the fibrations on the 
appropriate normal bundles.

\section{Compositions and Parametrices}




To prove closedness under composition we require to assume that a groupoid $\G$ which integrates the Lie structure on $M$ and a
groupoid $\G_{\partial}$ which integrates the Lie structure on $Y$ are chosen in the following sense.
Precisely, we construct for the given Lie structures on $M$ and $Y$ respectively integrating groupoids $\G$ and $\G_{\partial}$ as well as a morphism from $\G \to \G_{\partial}$ and a morphism from $\G_{\partial} \to \G$ 
in the category of Lie groupoids.
These morphisms are described using \emph{correspondences} of groupoids, see \cite{MO}.
\begin{Exa}
Consider the example of the algebroids $\A = TM, \ \A_{\partial} = TY$, i.e. the Lie structures consisting of \emph{all vector fields}.  
The pair groupoids $M \times M \rightrightarrows M, \ Y \times Y \rightrightarrows Y$ and also the \emph{path groupoids} (see \cite{LN}, example 2.9) $\P_M \rightrightarrows M, \ \P_Y \rightrightarrows Y$ integrate these algebroids. 
\label{Exa:corr}
\end{Exa}
For several Lie structures, groupoids and correspondences with good geometry exist, e.g. the Lie structure of $b$-vector fields or the structure of fibered cusp vector fields \cite{B}.
The following results hold for Lie structures of this type.

\begin{Thm}
The class of extended Boutet de Monvel operators $\B_{2\V}^{0,0}(M, Y)$ is closed under composition and adjoint, hence
$\B_{2\V}^{0,0}(M, Y)$ forms an associative $\ast$-algebra. 
\label{Thm:closed1} 
\end{Thm}

We define the class of \emph{truncated} Boutet de Monvel operators as follows.
The restriction $r^+$ to the interior $\mathring{X}_0 := X_0 \setminus Y_0$ and the extension by zero operator $e^+$ are given on the manifold level by
\[
\xymatrix{
L^2(M_0) \ar@/^1pc/[r]^{r^{+}} & \ar@/-0pc/[l]^-{e^{+}} L^2(\mathring{X}_0) 
} 
\]

with $r^{+} e^{+} = \id_{L^2(\mathring{X}_0)}$ and $e^{+} r^{+}$ being a projection onto a subspace of $L^2(M_0)$.
We define
\[
\End\begin{pmatrix} C_c^{\infty}(M_0) \\ \oplus \\ C_c^{\infty}(Y_0) \end{pmatrix} \supset \B_{2\V}^{m,0}(M, Y) \ni A = \begin{pmatrix} P + G & K \\ T & S \end{pmatrix} \mapsto \C(A) = \begin{pmatrix} r^{+} (P + G) e^{+} & r^{+} K \\ T e^{+} & S \end{pmatrix} \in \End\begin{pmatrix} C_c^{\infty}(X_0) \\ \oplus \\ C_c^{\infty}(Y_0) \end{pmatrix}.
\]

\begin{Def}
The class of \emph{truncated operators} is for $m \leq 0$ defined as $\B_{\V}^{m,0}(X, Y) := \C \circ \B_{2\V}^{m,0}(M, Y)$.
\label{Def:truncated}
\end{Def}

To show closedness under composition we use the longitudinally smooth structure of the integrating groupoids
as well as the previous Theorem.
This enables us to state the second main result.

\begin{Thm}
The calculus $\B_{\V}^{0,0}(X, Y)$ is closed under composition and adjoint. 
\end{Thm}

A priori, the inverse of an invertible Boutet de Monvel operator will not be contained in our calculus due to the definition via compactly supported distributional kernels. 
We define a completion $\overline{\B}_{\V}^{-\infty,0}(X, Y)$ of the residual Boutet de Monvel operators with regard to the family of norms
of operators $\L\left(\begin{matrix} H_{\V}^t(X) \\ \oplus \\ H_{\W}^t(Y) \end{matrix}, \ \begin{matrix} H_{\V}^r(X) \\ \oplus \\ H_{\W}^r(Y) \end{matrix}\right)$ on Sobolev spaces, cf. \cite{ALNV}.
Define the completed algebra of Boutet de Monvel operators as
\[
\overline{\B}_{\V}^{0,0}(X, Y) = \B_{\V}^{0,0}(X, Y) + \overline{\B}_{\V}^{-\infty, 0}(X, Y).
\]
The resulting algebra contains inverses and has favorable algebraic properties, e.g. it is spectrally invariant, \cite{B}. 
We obtain a parametrix construction after defining a notion of \emph{Shapiro-Lopatinski ellipticity}.
The indicial symbol $\R_{F}$ of an operator $A$ on $X$ is an operator $\R_{F}(A)$ defined as the restriction to a singular hyperface $F \subset X$ (see \cite{ALN}). 
Note that if $F$ intersects the boundary $Y$ non-trivially we obtain in this way a non-trivial Boutet de Monvel
operator $\R_F(A)$ defined on the Lie manifold $F$ with boundary $F \cap Y$. 

\begin{Def}
\emph{i)} We say that $A \in \overline{\B}_{\V}^{0,0}(X, Y)$ is \emph{$\V$-elliptic} if the principal symbol $\sigma(A)$ and the principal boundary symbol $\sigma_{\partial}(A)$ are 
both pointwise invertible.

\emph{ii)} A $\V$-elliptic operator $A$ is \emph{elliptic} if $\R_F(A)$ is pointwise invertible for each hyperface $F \subset X$. 
\label{Def:elliptic}
\end{Def}

\begin{Thm}
\emph{i)} Let $A \in \overline{\B}_{\V}^{0,0}(X, Y)$ be $\V$-elliptic. There is a parametrix $B \in \overline{\B}_{\V}^{0,0}(X, Y)$ of $A$, in the sense
\[
I - AB \in \overline{\B}_{\V}^{-\infty, 0}(X, Y), \ I - BA \in \overline{\B}_{\V}^{-\infty, 0}(X, Y). 
\]


\emph{ii)} Let $A \in \overline{\B}_{\V}^{0,0}(X, Y)$ be elliptic. There is a parametrix $B \in \overline{\B}_{\V}^{0,0}(X, Y)$ of $A$ up to 
compact operators
\[
I - AB \in \K\begin{pmatrix} L_{\V}^2(X) \\ \oplus \\ L_{\W}^2(Y) \end{pmatrix}, \ I - BA \in \K\begin{pmatrix} L_{\V}^2(X) \\ \oplus \\ L_{\W}^2(Y) \end{pmatrix}. 
\]

\label{Thm:parametrix}
\end{Thm}

{ \small {

}}

\end{document}